\documentstyle[12pt]{amsart}

\makeatletter                                            
\renewcommand{\@begintheorem}[2]{                        
\rm \trivlist \item [\hskip \labelsep {\bf #2\ \ #1.}]   
                                }                        
\makeatother                                             

\makeatletter

\DeclareFontFamily{U}{cyr}{}
\DeclareFontShape{U}{cyr}{m}{n}{
  <5> wncyr5 <6> wncyr6 <7> wncyr7 <8> wncyr8 <9> wncyr9 <10->
wncyr10}{}
\DeclareMathAlphabet{\mathcyr}{U}{cyr}{m}{n}

\input cyracc.def
\input epsf

\newcommand{\ts}{\vspace{\baselineskip}\noindent{\bf Proof.}$\;\;$}
\newcommand{\ZZ}{{\bf Z}}
\newcommand{\QQ}{{\bf Q}}
\newcommand{\RR}{{\bf R}}
\newcommand{\CC}{{\bf C}}
\newcommand{\FF}{{\bf F}}

\newcommand{\PP}{{\bf P}}

\newcommand{\cC}{{\cal C}}
\newcommand{\cE}{{\cal E}}
\newcommand{\cF}{{\cal F}}
\newcommand{\cG}{{\cal G}}

\newcommand{\cO}{{\cal O}}
\newcommand{\cR}{{\cal R}}

\newcommand{\et}{{\text{\'et}}}
\newcommand{\rmod}{{\mbox{$\;{\rm mod}\;$}}}

\newcommand{\TSh}{{\mathcyr{\cyracc Sh}}}
\newcommand{\pE}{{'\!E}}

\begin{document}

\title{Some remarks on Brauer groups of K3 surfaces}
\author{Bert van Geemen}
\address{Dipartimento di Matematica, Universit\`a di Milano,
  via Saldini 50, I-20133 Milano, Italia}
 \email{geemen@mat.unimi.it}

\begin{abstract}
We discuss the geometry of the genus one fibrations associated to an
elliptic fibration on a K3 surface.
We show that the two-torsion subgroup of
the Brauer group of a general elliptic fibration
is naturally isomorphic to the two torsion of the Jacobian
of a curve associated to the fibration.
We remark that this is  related to Recillas' trigonal construction.
Finally we discuss the two-torsion in the Brauer group
of a general K3 surface with a polarisation of degree two.
\end{abstract}

\maketitle

Elements of the Brauer group of a projective variety have
various geometric incarnations, like Azumaya algebras and Severi-Brauer
varieties. Moreover,
an element $\alpha$ in the Brauer group of a projective
K3 surface $X$ determines
a polarized Hodge substructure $T_{<\alpha>}$
of the transcendental lattice $T_X$ of $X$.
We are interested in finding geometric
realizations of the Hodge structure $T_{<\alpha>}$ and in relating
such a realization to the other incarnations.
In this paper we collect some examples and add a few new results.
Recently, there has been much interest in Brauer groups of K3 surfaces
(\cite{Cnfm}, \cite{DP}, \cite{W}), but our aim is much more modest.

For certain elements $\alpha\in Br(X)$ the Hodge structure
$T_{<\alpha>}$ is the transcendental lattice of another K3 surface $Y$.
The surface $Y$ is not necessarily unique.
Two K3 surfaces with
isomorphic transcendental lattice are said to be Fourier-Mukai partners.
Results from Mukai and Orlov show that
if $Z$ is a K3 surface then $T_Y$ and $T_Z$ are Hodge isometric if and
only if
the bounded derived categories of coherent sheaves on $Y$ and $Z$ are
isomorphic. The number of Fourier-Mukai partners of
a given K3 surface is determined in \cite{HLOY}, see also \cite{St}.

In sections \ref{brgr} and \ref{brk3} we briefly recall some basic
definitions.
In section \ref{brelk3} we consider the case when the K3 surface
$X$ has an elliptic fibration (with a section).
In that case $T_{<\alpha>}$ is always the transcendental lattice
of a K3 surface $Y$, but $Y$ is far from unique in general.
These surfaces are studied in section \ref{g1f}.
They have a fibration in genus one curves (without a section).
In fact, any pair of elements $\pm\alpha\in Br(X)$
determines a unique fibration in genus one curves on such a surface.
Somewhat surprisingly, $Y$ often has more than one genus one fibration.
We use well-known results on K3 surfaces to determine
the number of these surfaces,
their automorphisms and the genus one fibrations.
In Example \ref{example}
we work out the numerology for certain cyclic subgroups of $Br(X)$.
In section \ref{exaY} we consider projective models for some genus one
fibrations.

In section \ref{brdo} we relate the 2-torsion in the Brauer group of a
double cover of a rational surface to the 2-torsion of the Jacobian of the
ramification curve, the main result, Theorem \ref{thmn},
is due to the referee. It greatly
generalizes a result in a previous version of this paper.
In section \ref{ttbr} we use Theorem \ref{thmn} to show that
the  2-torsion in the Brauer group
of a general elliptic fibration $X\rightarrow \PP^1$
is naturally
isomorphic to the
2-torsion in the Jacobian of a certain curve in $X$.
We observe that the pairs $X$, $Y$,
where $Y$ is the surface associated to an element $\alpha$ of order two
in $Br(X)$, are closely related to
Recillas' trigonal construction. This construction gives a bijection between
trigonal curves with a point of order two and
tetragonal curves.
In section \ref{cobu} we remark that classical formulas allow one to
give the conic bundle, and hence the Azumaya algebra,
on $X$ associated to $\alpha$.

In the last section we consider K3 surfaces with N\'eron-Severi group
of rank one and we characterize the elements $\alpha$
of order two in the Brauer group such that
$T_{<\alpha>}$ is the transcendental lattice of another K3 surface.
In case the K3 surface is a double plane, we also consider the other
elements of order two.
It is well-known that some of these are related to cubic fourfolds,
in this way one obtains a geometric realization of $T_{<\alpha>}$.
We briefly discuss the remaining case in \ref{remaining}.

\section{Brauer groups}\label{brgr}

\subsection{Basic definitions and properties}
The main references for
Brauer groups of algebraic varieties and, more generally, schemes are
Grothendieck's papers \cite{G} and Milne's book \cite{M}, useful
summaries are given in \cite{Ct}, \cite{Cnfm}, \cite{DP}.
For recent work on non-projective complex varieties see \cite{HS}.

The cohomological Brauer
group of a scheme $X$ is $Br'(X)_{\et}=H^2_{\et}(X,\cO_X^*)$.
For a smooth complex algebraic variety $X$ one has
$Br'(X)_{\et}\cong H^2(X,\cO_X^*)_{tors}$, the subgroup of torsion
elements in $H^2(X,\cO_X^*)$.

 An Azumaya algebra is a
sheaf of $\cO_X$-algebras over $X$,
which is locally isomorphic to $\cE nd(\cE)$ for a locally free sheaf
$\cE$.
Two Azumaya algebras are (Morita) equivalent if
${\cal A}\otimes {\cE}nd(\cE)\cong{\cal A'}\otimes
{\cE} nd({\cE}')$
where ${\cE},{\cE}'$ are locally free sheaves on $X$.
The Brauer group $Br(X)$ of a complex algebraic variety
$X$ is the group of equivalence classes of Azumaya algebras
with tensor product.

The Brauer group also classifies Severi-Brauer varieties, these are
projective bundles $P$ over $X$ (\cite{Artin}).

For smooth curves one has $Br(X)=Br'(X)_{\et}=0$. For a smooth surface
one has $Br(X)=Br'(X)_{\et}$ and recently A.J.\ de Jong \cite{dJ} proved that
an element of order $n$ in $Br(X)$ has a representative which is an
Azumaya algebra locally isomorphic to ${\cE}nd(\cE)$
for a locally free sheaf $\cE$ of rank $n$.

\subsection{Brauer groups of elliptic fibrations}\label{oggsha}
Let $X$ be a smooth algebraic surface
with an elliptic fibration $p:X\rightarrow \PP^1$ (with
a section). The Brauer group of $X$ is isomorphic to the
Tate-Shafarevich group
$\TSh(X)=H^1_{\et}(\PP^1,X^\sharp)$
where $X^\sharp$ is the sheaf of groups on
$\PP^1$ of local sections of $p_X:X\rightarrow \PP^1$
(\cite{CD} Ch.\ 5, \cite{DG}, \cite{Ct}, Theorem 4.4.1,
cf.\ \cite{FM}, 1.5.1):
$$
Br(X)=\TSh(X).
$$
A non-trivial element of order $n$ in $\TSh(X)$ is a genus one fibration
$Y\rightarrow \PP^1$ (without a section) and $n$ is
the minimal degree of a multisection of $p:Y\rightarrow \PP^1$.
There is a birational action of $X$ on $Y$.
The genus one fibrations corresponding to
$\alpha$ and $-\alpha\in \TSh(X)$
are isomorphic as fibrations
(\cite{Ct}, Remark 4.5.3), but the action of $X$ on them
is different. For the spaces corresponding to $k\alpha$, $k\in\ZZ$, see
\cite{DP}, \S 2.2, \cite{Ct}, Thm 4.5.2.
The fibration $Y/\PP^1$ has a relative Jacobian fibration which
is $X/\PP^1$.

\section{The Brauer group of a K3 surface}\label{brk3}

\subsection{}\label{torB}
For a smooth surface with $H^3(X,\ZZ)=0$, the exponential sequence shows
that
$$
Br(X)=H^2(X,\cO_X^*)_{tors}\cong
\left(H^2(X,\ZZ)/c_1(Pic(X))\right)\otimes(\QQ/\ZZ),
$$
where $c_1:Pic(X)\rightarrow H^2(X,\ZZ)$ is the first Chern class.
In case $X$ is K3 surface, $Pic(X)\cong c_1(Pic(X))=NS(X)$, the
N\'eron-Severi group of $X$. The transcendental lattice of $X$ is
the sublattice of $H^2(X,\ZZ)$ perpendicular to $NS(X)$ w.r.t.\ the
intersection pairing:
$$
T_X=NS(X)^\perp=\{v\in H^2(X,\ZZ):vn=0\;\mbox{for all}\; n\in NS(X)\,\}.
$$
The intersection pairing is unimodular,
in fact,
$$
H^2(X,\ZZ)\cong \Lambda_{K3}:=U\oplus U\oplus U\oplus E_8(-1)\oplus E_8(-1)
$$
where $U$ is the hyperbolic plane
(cf.\ \cite{BPV}, VIII.1).
In particular $H^2(X,\ZZ)$ is self-dual, and
as $H^2(X,\ZZ)/T_X$ is torsion free
(i.e.\ $T_X$ is a primitive sublattice)
we get
$$
H^2(X,\ZZ)/NS(X)\stackrel{\cong}{\longrightarrow}T_X^*:=Hom(T_X,\ZZ),
\qquad v+NS(X)\longmapsto [t\mapsto vt]
$$
for $v\in H^2(X,\ZZ)$, $t\in T_X$. Thus for a K3 surface $X$ we obtain:
$$
Br(X)\cong T_X^*\otimes(\QQ/\ZZ)\cong Hom(T_X,\QQ/\ZZ).
$$

\subsection{Hodge structures}\label{hodgek3}
Let $X$ be a K3 surface.
The Hodge structure on $H^2(X,\ZZ)$ is determined by the
one dimensional subspace
$H^{2,0}(X)=<\omega_X > \subset H^2(X,\ZZ)\otimes \CC$.
A Hodge isometry $f:H^2(X,\ZZ)\rightarrow H^2(Y,\ZZ)$
is a $\ZZ$-linear map which preserves the intersection forms and
whose $\CC$-linear extension maps $H^{2,0}(X)$ to $H^{2,0}(Y)$.

The (weak) Torelli theorem for projective K3 surfaces
(\cite{BPV}, VIII.6.3)
asserts that if  $H^2(X,\ZZ)$ and $H^2(Y,\ZZ)$ are Hodge isometric,
where $X$ and $Y$ are K3 surfaces, then $X$ and $Y$
are isomorphic.

The period $\omega_X$ of $X$ satisfies $\omega_X\omega_X=0$ and
$\omega_X\bar{\omega}_X>0$.
The `surjectivity of the period map' (\cite{BPV}, VIII.14) asserts
that given an $\omega\in \Lambda_{K3}\otimes\CC$ with
$\omega\omega=0$ and $\omega\bar{\omega}>0$ there exists a K3 surface $X$
and a Hodge isometry
$f:(\Lambda_{K3},\omega)\rightarrow H^2(X,\ZZ)$
(so $<f_\CC(\omega)>=H^{2,0}(X)$).
The N\'eron-Severi group of $X$ is isomorphic to the lattice
$\{v\in \Lambda_{K3}:v\omega=0\}$.
The K3 surface $X$ is projective iff $NS(X)$ contains a class $v$
with $v^2>0$, equivalently, the signature of $T_X=NS(X)^\perp$
is $(2+,n-)$ for some $n$, $0\leq n\leq 19$.

The transcendental lattice $T_X$ is a polarized Hodge structure
(note that $<\omega_X>\subset T_X\otimes\CC$).
A sublattice $\Gamma\subset T_X$ of finite index
is also a polarized Hodge structure, in fact
$<\omega_X>\subset \Gamma\otimes\CC=T_X\otimes\CC$.
In case $\Gamma$
admits a primitive embedding $i:\Gamma\hookrightarrow \Lambda_{K3}$,
the vector $i_\CC(\omega_X)\subset \Lambda_{K3}\otimes\CC$
defines a K3 surface $Y=Y(\Gamma,i)$ with period
$\omega_Y=i_\CC(\omega_X)$ and $T_Y= i(\Gamma)$.
If $X$ is projective, $Y$ is also projective (the transcendental lattices
have the same signatures).

\subsection{Moduli spaces}
Conversely, let $Y$ be a K3 surface and let $T_Y \hookrightarrow T$ be a
Hodge isometry such that $T/T_Y$ is cyclic.
Mukai \cite{Mmb} has shown that
there is a K3 surface $X$,
actually $X$ is a moduli space of sheaves on $Y$, such that $T\cong T_X$.

Kond\=o \cite{K1}, \cite{K2} studied the map defined by such inclusions
from the moduli space of K3 surfaces with a
polarization of degree $dn^2$ to the one of degree $d$,
this map has a high degree. We are interested in going in the opposite
direction, part of the extra data needed for the inverse is a cyclic subgroup
of the Brauer group of $X$.

\subsection{Brauer groups and Hodge substructures}\label{brhod}
An element $\alpha$ of order $n$ in $Br(X)$ defines a surjective
homomorphism $\alpha:T_X\rightarrow\ZZ/n\ZZ$ (cf.\ \ref{torB})
and thus defines a sublattice of index $n$ of $T_X$.
Conversely, any sublattice $\Gamma'$
of index $n$
in $T_X$ with cyclic quotient group $T_X/\Gamma'$
determines a subgroup
of order $n$ in $Br(X)$.
The sublattice defined by the cyclic subgroup $\langle\alpha\rangle$
of $Br(X)$ generated by $\alpha$ is denoted by
$$
T_{<\alpha>}=\ker(\alpha:T_X\longrightarrow \QQ/\ZZ)\qquad
(\alpha\in Br(X)=Hom(T_X,\QQ/\ZZ)).
$$

In case there exists a primitive embedding
$i:T_{<\alpha>}\hookrightarrow \Lambda_{K3}$,
we obtain a K3 surface $Y=Y(T_{<\alpha>},i)$ as in \ref{hodgek3}.
Such an embedding $i:T_{<\alpha>}\hookrightarrow \Lambda_{K3}$
is not unique up to isometry in general.
In fact, even the orthogonal complement of the images
under the various embeddings need not be isomorphic as lattices
(cf.\ \cite{O}, \cite{HLOY}, \cite{St}).

\subsection{Lattices}\label{lattices}
We recall some basic definitions (\cite{Dol}, \cite{N}).
Given a lattice $(L,<,>)$, that is, a free $\ZZ$-module $L$
with a symmetric
bilinear form $<,>:L\times L\rightarrow\ZZ$, the dual lattice $L^*$
is defined as
$$
L^*=Hom(L,\ZZ)=\{x\in L\otimes\QQ:\;<x,y>\in\ZZ\;\mbox{for all}\;y\in L\}.
$$
The discriminant group of $L$ is $L^*/L$ and the discriminant form
of an even lattice (so $<x,x>\in2\ZZ$ for all $x\in L$) is the map
$$
q:L^*/L\longrightarrow \QQ/2\ZZ,\qquad q(x+L)=<x,x>\;\rmod 2\ZZ
$$
where the bilinear form is extended $\QQ$-bilinearly to $L^*$.

\section{Brauer groups of elliptic fibrations on K3 surfaces}\label{brelk3}

\subsection{The transcendental lattice of an elliptic K3 fibration}
A general elliptic K3 surface $f:X\rightarrow \PP^1$ with a section
has N\'eron-Severi group isometric to $U$, the hyperbolic plane.
It is spanned by the class of a fiber
$f$ and the section $s$. As $f^2=0,fs=1,s^2=-2$, we get
$f(f+s)=1,(f+s)^2=0$ hence $f$ and $f+s$ are a standard basis of $U$.
The elliptic fibration, and the section, on $X$ are unique
(cf.\ \ref{1,0}).
We will simply write $X/\PP^1$ for this unique Jacobian fibration.

There is, up to isometry, a unique embedding
$U\hookrightarrow \Lambda_{K3}$ (\cite{N}, Thm 1.14.4) hence
the transcendental lattice $T_X=NS(X)^\perp$ of $X$ is isometric to:
$$
T_X\cong \Gamma:= U\oplus U\oplus E_8(-1)\oplus E_8(-1).
$$
In particular, $\Gamma$ is unimodular and even.

\subsection{The Brauer group of a general elliptic K3 fibration}
An element $\alpha$ of order $b$ in $Br(X)$ defines a surjective
homomorphism $\alpha:\Gamma\rightarrow\ZZ/b\ZZ$ (cf.\ \ref{torB})
and determines the sublattice
$T_{<\alpha>}=\ker(\alpha:\Gamma \longrightarrow \QQ/\ZZ)$
of index $b$ in $\Gamma$.
As a first step
towards understanding the Brauer group we have the following
classification.

\subsection{Proposition.}\label{lemed}
Let $\Gamma'\subset\Gamma$ be a sublattice of index $b$ such that
$\Gamma/\Gamma'$ is cyclic.
Then there exists a $c\in\ZZ$ such that
$$
\Gamma'\cong\Gamma_{b,c}:= \Lambda_{b,c}\oplus U\oplus E_8(-1)^2,
$$
where $\Lambda_{b,c}$ is the indefinite lattice of rank two defined as:
$$
\Lambda_{b,c}:=\left\{\ZZ^2,
(0,b,2c):=\;\left(\begin{array}{cc} 0&b\\b&2c\end{array}\right)\right\}.
$$
The lattice $\Gamma_{b,c}$ is determined, up to isometry,
by its rank, signature and discriminant form.

Moreover, if $gcd(b,2c)=1$ then:

\noindent
$\Gamma_{b,c}\cong \Gamma_{b,d}$ iff $c\equiv d\rmod b$ and
there is an $a\in\ZZ$ such that $d\equiv a^2c\rmod b^2$.

\ts
Since $\Gamma$ is unimodular, there is a $\gamma\in \Gamma$
such that
$$
\Gamma'=\ker\Big(\Gamma\longrightarrow\ZZ/b\ZZ,\qquad
x\longmapsto \langle x,\gamma\rangle\;\rmod\,b\ZZ\Big),
$$
and we may assume that $\gamma$ is a primitive vector.

Let $\gamma^2=-2c$ and assume that $c\neq 0$.
As $\Gamma$ is even, unimodular, of signature $(2,18)$, the embedding of
the sublattice $M=\ZZ\gamma$ into $\Gamma$ is unique up to isomorphism
(\cite{N}, Thm 1.14.4). In particular, we may assume,
up to isometries of $\Gamma$, that
$\gamma=(1,-c)\in U\subset\Gamma$,
where $U$ is the first summand of $\Gamma$. Then $\Gamma'$ is generated
by $(0,b),(1,c)\in U\subset\Gamma$
and the unimodular summand $U\oplus E_8(-1)^2$.
Thus $\Gamma'\cong \Gamma_{b,c}$.

By \cite{N}, Cor.\ 1.13.3,
the lattice $\Gamma_{b,c}$ is determined up to isometry by its rank,
signature and discriminant form. The discriminant group $A_{b,c}$
of $\Gamma_{b,c}$ is the one of $\Lambda_{b,c}$, which has order
$disc(\Lambda_{b,c})=b^2$.
If $gcd(b,2c)=1$,  $A_{b,c}$ is cyclic of order $b^2$ with generator
$(-2c,b)/b^2$ and the discriminant form $q$ takes the value
$q((-2c,b)/b^2)=-2c/b^2\in\QQ/2\ZZ$. If $\Gamma_{b,d}$ has isomorphic
discriminant group, then $gcd(b,2d)=1$ (else $A_{b,d}$ is not cyclic)
and an isomorphism $A_{b,c}\rightarrow A_{b,d}$ maps the generator
$(-2c,b)/b^2$ to the generator $a(-2d,b)/b^2$, for some integer $a$
with $gcd(a,b^2)=1$.
The isomorphism should respect discriminant forms, hence
$-2c/b^2=a^2(-2d/b^2)$ in $\QQ/2\ZZ$, so $c\equiv a^2d\rmod b^2$.
Viceversa, such an $a\in\ZZ$ defines an isomorphism $A_{b,c}\cong A_{b,d}$.

In case $\gamma^2=0$, choose any $\delta'\in \Gamma$ with
$\gamma\delta'=1$, such a $\delta'$ exists since $\gamma$ is primitive
and $\Gamma$ is unimodular.
If $(\delta')^2=2d$, define $\delta:=\delta'-d\gamma$, then
the sublattice $H:=\ZZ\gamma+\ZZ\delta$
is isomorphic to $U$. Hence $\Gamma\cong H\oplus H^\perp$ where
$H^\perp$ is an even unimodular lattice
(isomorphic to $U\oplus E_8(-1)^2$).
Then $\Gamma'$ is generated by $\gamma$, $b\delta$ and $H^\perp$,
so $\Gamma'\cong \Gamma_{b,0}$.
\qed

\subsection{Corollary}
Let $b>2$ be a prime number. Then there are three isomorphism
classes of lattices $\Gamma_{b,c}$:
those with $c\equiv 0\rmod b$, those where $c$
is congruent to a non-zero square mod $b$ and those where $c$ is not
congruent to a square mod $b$.

\ts
The condition $c\equiv a^2d\rmod b^2$ implies of course
$c\equiv a^2d\rmod b$, and it is easy to see that
the converse holds as well.
\qed

\subsection{The case $b=2$}\label{b=2}
In case $b=2$ there are two isomorphism classes with representatives
$\Gamma_{2,0}$ and $\Gamma_{2,1}$ respectively.
As in the proof of Proposition \ref{lemed}
let $\alpha\in Br(X)_2=Hom(\Gamma,\ZZ/2\ZZ)$,
correspond to the homomorphism $x\mapsto <x,\gamma>$
defined by $\gamma\in \Gamma$.
The quadratic form on $\Gamma$ induces a non-degenerate quadratic form:
$$
\Gamma/2\Gamma \cong \FF_2^{20}\longrightarrow \FF_2,\qquad
\bar{\gamma}\longmapsto \langle\gamma,\gamma\rangle/2 \,\rmod 2\ZZ.
$$
This is an `even' quadratic form.
In particular, it has $2^{9}(2^{10}+1)$
zeroes in $\FF_2^{20}$, including $\bar{\gamma}=0$. Thus
there are $2^{9}(2^{10}+1)-1$ elements in $Br_2(X)$ with
$T_{<\alpha>}\cong \Gamma_{2,0}$
and there are $2^{9}(2^{10}-1)$
elements with $T_{<\alpha>}\cong \Gamma_{2,1}$.

\subsection{Embeddings of $T_{<\alpha>}$ in $\Lambda_{K3}$}
Given an $\alpha\in Br(X)$ of order $n$, the sublattice $T_{<\alpha>}$
of $T_X$ is isomorphic to some $\Gamma_{b,c}$ as in Proposition
\ref{lemed}.
Any lattice $\Gamma_{b,c}$ has a primitive embedding into
the K3 lattice. For example, one can simply map the standard basis
of $\Lambda_{b,c}$ to $((1,0),(b,0))$, $((0,0),(2c,1))\in U^2$,
the first two copies of $U$ in $\Lambda_{K3}$, and extend it by
the identity on $U\oplus E_8(-1)^2$ to all of $\Gamma_{b,c}$.
In particular, there is a K3 surface $Y$ with $T_Y\cong T_{<\alpha>}$
(cf.\ \ref{brhod}), which is not surprising in view of the
Ogg-Shafarevich theory (cf.\ \ref{oggsha}).

The orthogonal complement of $i(\Gamma_{b,c})$ depends
in general on the primitive embedding $i$.
In fact, we will see that it is isomorphic to $\Lambda_{b,-a^2c}$
for some $a\in\ZZ$
but the classification of the lattices $\Lambda_{b,d}$
up to isomorphism is different
from the one of the $\Gamma_{b,c}$ (cf.\  Proposition \ref{lemed}).
Indeed, the summand $U\oplus E_8(-1)^2$ in $\Gamma_{b,c}$ leads to fewer
isomorphism classes. Therefore it can happen that
$\Gamma_{b,c}\cong \Gamma_{b,d}$ but
$\Lambda_{b,c}\not\cong \Lambda_{b,d}$.

\subsection{Proposition}\label{Lbc}
Let $b$ be a non-zero integer and let
$i: \Gamma_{b,c}\hookrightarrow \Lambda_{K3}$
be a primitive embedding. Then there is a non-zero integer $a$ such that
$$
i(\Gamma_{b,c})^\perp \cong \Lambda_{b,-a^2c}.
$$
Conversely, for any non-zero integer $a$ there is a primitive embedding
such that $i(\Gamma_{b,c})^\perp \cong \Lambda_{b,-a^2c}$.

For $b,c,d,n\in\ZZ$ with $b\neq 0, cd\equiv 1\rmod b$
we have:
$$
\Lambda_{b,c}\cong \Lambda_{-b,c}\cong \Lambda_{b,c-nb}\cong \Lambda_{b,d}.
$$
Moreover, if $gcd(b,c)=1$ and $0\leq c<d\leq b-1$ then:

$\Lambda_{b,c}\cong \Lambda_{b,d}\;$ iff $\;cd\equiv 1 \rmod b$.

\ts
As $\Lambda_{K3}$ is unimodular, the lattice $i(\Gamma_{b,c})^\perp$
has the same discriminant $-b^2$ as $\Gamma_{b,c}$
and it has the opposite discriminant form.
Moreover it has rank two and it is indefinite. Writing the quadratic form
on $i(\Gamma_{b,c})^\perp$ as $kx^2+2lxy+my^2$, its discriminant is
$b^2=l^2-km$, hence
there is an isotropic vector $v$ in $i(\Gamma_{b,c})^\perp$,
i.e.\ $v\neq 0, v^2=0$.
We can choose $v$ to be primitive. Thus there is a
$w\in i(\Gamma_{b,c})^\perp$ such
that $v,w$ is a basis of $i(\Gamma_{b,c})^\perp$
and the quadratic form on this basis has $k=0$. Then $l=\pm b$,
hence
$i(\Gamma_{b,c})^\perp \cong \Lambda_{b,d}$ for some $d$.
The discriminant form of $i(\Gamma_{b,c})^\perp$ is the opposite of the
one of $\Gamma_{b,c}$ which we determined in the proof of Proposition
\ref{lemed} and we get $i(\Gamma_{b,c})^\perp \cong \Lambda_{b,-a^2c}$.
For the converse, use \cite{N} Prop.\ 1.4.1,
which shows that there exists an even unimodular lattice
containing $\Gamma_{b,c}\oplus\Lambda_{b,-a^2c}$ as sublattice of
finite index. The classification of even unimodular lattices implies
that this overlattice is isomorphic to $\Lambda_{K3}$
and this gives the desired embedding $i$.

Let $e_1,e_2$ be the standard basis of $\Lambda_{b,c}$, so $e_1^2=0$,
$e_1e_2=b$ and $e_2^2=2c$.
Let $f_1:=-e_1,f_2=e_2$, then $f_i^2=e_i^2$ and $f_1f_2=-b$,
so $\Lambda_{b,c}\cong \Lambda_{-b,c}$.
Next, let $f_1=e_1, f_2:=-ne_1+e_2$ then
$f_1^2=0$, $f_1f_2=b$, $f_2^2=2(c-nb)$, hence
$\Lambda_{b,c}\cong \Lambda_{b,c-nb}$.
Moreover, if $cd\equiv 1\rmod b$, then
$cd-be=1$ for an $e\in\ZZ$ and we have:
$$
^tA\left(\begin{array}{cc} 0&b\\b&2c\end{array}\right) A=
\left(\begin{array}{cc} 0&b\\b&2d\end{array}\right),\qquad
{\rm with}\qquad
A=A_{b,c,d}=\left(\begin{array}{cc} -c&-e\\b&d\end{array}\right)
\in GL(2,\ZZ).
$$
The last statement follows by explicitly considering two by two integral
invertible matrices $A$ satisfying $^tA(0,b,2c)=(0,b,2d)A^{-1}$.
(cf.\ \cite{St}, Lemma 4.1).
\qed

\section{Genus one fibrations on K3 surfaces}\label{g1f}

\subsection{} In the previous section
we considered a general elliptic fibration $X\rightarrow \PP^1$
on a K3 surface.
We found that for $\alpha\in Br(X)$ of order $b$ the lattice
$T_{<\alpha>}$ is isomorphic to the transcendental lattice
of a K3 surface $Y$ and that the N\'eron-Severi group of $Y$
is isomorphic to $\Lambda_{b,c}$ (as in Proposition \ref{lemed}),
but $c$ is in general not uniquely determined by $T_{<\alpha>}$.
We may, and will, assume that $b>0$ and $0\leq c <b$
from now on.

From Ogg-Shafarevich theory we know that $Y$ has a genus one fibration
$p_Y:Y\rightarrow \PP^1$ with a multisection of degree $b$ and the
relative Jacobian fibration is $X\rightarrow \PP^1$.
We will now see that in general $Y$ has two genus one fibrations
and that the relative Jacobians of these fibrations are not always
isomorphic.

\subsection{Remark} If $Z\rightarrow \PP^1$ is a
general genus one fibration on an  algebraic K3 surface,
then it is easy to see that $NS(Z)\cong\Lambda_{b,c}$ for some $b>0$.

\subsection{The K\"ahler cone of $Y_{b,c}$}\label{ybc}
To study the geometry (genus one fibrations, automorphisms and
projective models) of a K3 surface $Y=Y_{b,c}$ with
$NS(Y)\cong \Lambda_{b,c}$
we need to know its K\"ahler cone
$C^+_{b,c}\subset NS(Y_{b,c})\otimes\RR$.
For this we need to determine the set
$\Delta\subset NS(Y)$ of (classes of) $(-2)$-curves:
$$
\Delta=\{D\in NS(Y):\; D>0,\;D\;{\rm irreducible},\;D^2=-2\}
$$
(\cite{BPV}, Chapter VIII, Cor.\ 3.8). The following easy lemma
supplies this information.

\subsection{Lemma} \label{kc} Let $Y_{b,c}$ be a K3 surface with
N\'eron-Severi group isomorphic to $\Lambda_{b,c}$, $0\leq c< b$.

Then there is an isomorphism of lattices $NS(Y_{b,c})\otimes\RR\cong
\Lambda_{b,c}\otimes\RR\cong \RR^2$
such that:
\begin{enumerate}
\item
if $c\neq b-1$,
$$
C^+_{b,c}=\{(x,y)\in \RR^2:\; y>0\;{\rm and}\;bx+cy>0\}\qquad (0\leq c<b-1),
$$
and there are no $(-2)$-curves on $Y_{b,c}$.
\item
if $c=b-1$,
$$
C^+_{b,b-1}=\{v=(x,y)\in \RR^2: v\in C\;{\rm and}\; bx+(b-2)y>0\}\qquad
{\rm and}\qquad N=(-1,1),
$$
where $N$ is the unique $(-2)$-curve on $Y_{b,b-1}$.
\end{enumerate}

\ts
Let $C\subset NS(Y)\otimes\RR$ be the connected component
of the set $\{v\in NS(Y)\otimes\RR:\;v^2>0\}$ which contains
an ample divisor. Any isometry $NS(Y_{b,c})\cong \Lambda_{b,c}$
maps $C$ to one of the two connected components of
the set $\{(x,y)\in\RR^2:y(bx+cy)>0\}$.
Since $-I$ is an isometry of
$\Lambda_{b,c}$, we can choose an isometry
$NS(Y)\cong \Lambda_{b,c}$ which identifies
$$
C=\{(x,y)\in \RR^2:\; y>0\;{\rm and}\;bx+cy>0\}.
$$

To determine $\Delta$, note that
$v^2=2y(bx+cy)$ and $0\leq c<b$ imply: $v^2=-2$ iff $c=b-1$
and $v=\pm (1,-1)$. Moreover, as $v^2=-2$, either $v$ or $-v$
is the class of an effective divisor and this must be the class of a
$(-2)$-curve. As the reflection
$s_v$ (defined by $s_v(w)=w+(w,v)v$)
is an isometry and $s_v(v)=-v$, we may assume that
$$
\Delta=\{(-1,1)\}\;\subset\;\Lambda_{b,b-1}.
$$
The K\"ahler cone $C^+$ of $Y$ is
$$
C^+=
\{v=(x,y)\in \RR^2: v\in C\;{\rm and}\; v n>0\;
\mbox{for all}\; n\in\Delta\}.
$$
(\cite{BPV}, Chapter VIII, Cor.\ 3.8). \qed

\subsection{Automorphisms} The Torelli theorem for K3 surfaces
(\cite{BPV}, VIII.11)
allows one to determine the automorphism group of a K3 surface.
It coincides (via the induced action on $H^2(X,\ZZ)$) with
the group of effective Hodge isometries, i.e.\ Hodge isometries
which map the K\"ahler cone into itself.
We will assume in the remainder of this section that the only
Hodge isometries of $T_X$ are $\pm I$. This is the case for a general
K3 surface with a transcendental lattice of rank at least 3, and thus
holds for the general $Y_{b,c}$.

\subsection{Lemma}\label{autL}\label{autK3}
Let $d=gcd(b,c)$
and $b'=b/d$, $c'=c/d$. Then the orthogonal group $O(\Lambda_{b,c})$
of the lattice $\Lambda_{b,c}$ is:
$$
O(\Lambda_{b,c})\cong\left\{
\begin{array}{rl}
\{\pm I,\pm J\} &{\rm if}\; c=0\;{\rm or}\;(c')^2 \equiv 1 {\rmod} b',\\
\{\pm I\}&{\rm else,}\\
\end{array}
\right.
$$
where $J$ permutes the two isotropic sublattices of $\Lambda_{b,c}$.

The general K3 surface $Y_{b,c}$ with N\'eron-Severi group isomorphic to
$\Lambda_{b,c}$ has trivial automorphism group except when:
$$
(b,c)=(1,0),\;(2,0),\quad{\rm or}\quad  b>2,\,c=1,
$$
in which case $Aut(Y_{b,c})\cong \ZZ/2\ZZ$.

\ts
By considering the action of $O(\Lambda_{b,c})$ on the primitive
isotropic vectors it is clear that this group has at most 4 elements.
In case $c^2\equiv 1\rmod b$ the map $A_{b,c,c}$
in the proof of Proposition \ref{Lbc}
gives an isometry,
of order two, of the lattice $\Lambda_{b,c}$.

Let $Y=Y_{b,c}$, let $T$ be its transcendental lattice let
$\Lambda=\Lambda_{b,c}\cong NS(Y)$ and let
$\Lambda^*/\Lambda\cong T^*/T$ be the discriminant group.
Any effective Hodge isometry of $H^2(Y,\ZZ)$ must preserve
the K\"ahler cone and, by assumption, it induces $\pm I$ on $T$.
Note that $-I\in Aut(NS(Y))$ does not preserve the K\"ahler cone.

In case $(b,c)=(1,0)$, $NS(Y)$ is the hyperbolic plane.
The subgroup of $Aut(NS(Y))$ preserving the K\"ahler cone is
trivial.
As $H^2(Y,\ZZ)=NS(Y)\oplus T$,
there is an effective Hodge isometry on $H^2(Y,\ZZ)$ which
induces the identity on the
$NS(Y)$ and minus the identity on $T$ and this is the only
non-trivial effective Hodge isometry.
The corresponding automorphism
is the inversion map on the elliptic fibration on $Y=Y_{1,0}$
(cf.\ \ref{1,0}).

In case $b>1$ and $c=0$,
let $\{I,J\}$ be
the subgroup of $Aut(NS(Y))$ preserving the K\"ahler cone.
Then $J$ permutes the basis vectors of
$\Lambda^*/\Lambda\cong (\ZZ/b\ZZ)^2$, whereas $-1_T$ acts as $-1$
on this group. Hence only in case $b=2$ there is a non-trivial
automorphism on $Y_{b,0}$: it induces the identity on the
N\'eron-Severi group and minus the identity on the transcendental lattice.
It is the covering involution for the 2:1 map
$Y_{2,0}\rightarrow \PP^1\times\PP^1$ in \ref{2,0}.

The case $0<c<b$ and $(c')^2\equiv 1\rmod b'$.
The automorphisms $-I,-A_{b',c',c'}$ of $NS(Y)=\Lambda$
do not preserve the K\"ahler cone.
The vectors $(1/b)(1,0)$, $(1/b)(-c',b')\in \Lambda\otimes\QQ$ are in
$\Lambda^*$ and $A_{b',c',c'}$ permutes them.
If $A_{b',c',c'}$ lifts to an automorphism of $H^2(Y,\ZZ)$
which is the identity on
$T$, we must have $(1/b)(1,0)\equiv (1/b)(-c',b')\rmod\ZZ^2$, hence
$b'=b$, so $d=1$ and also $c=c'$. Moreover, as $0\leq c<b$
we get $c=b-1$.
In this case $A_{b,c,c}$ does not preserve the
K\"ahler cone.
If $A_{b',c',c'}$ lifts to an automorphism of $H^2(Y,\ZZ)$
which is minus the identity on $T$,
we must have $-(1/b)(1,0)\equiv (1/b)(-c',b')\rmod\ZZ^2$,
hence $b'=b,c'=c$ as before and $c'=1$.
In this case $Y$ has indeed an automorphism of order two which is
the covering involution for the 2:1 map $\phi_h:Y\rightarrow \PP^2$
given by $h=(0,1) \in C^+$ with $h^2=2$,
cf.\ the examples \ref{2,1}, \ref{3,1} below.
\qed

\subsection{The elliptic fibrations on $Y_{b,c}$}\label{ellfibc}
The elliptic fibrations on $Y_{b,c}$ correspond
to the primitive isotropic classes in the
closure of the K\"ahler cone. Thus there are two such fibrations
if $c\neq b-1$:
$$
\epsilon_1=(1,0),\quad \epsilon_2=(-c',b')\quad \in C^+_{b,c},
\qquad \epsilon_1^2=\epsilon_2^2=0,
$$
with $b'=b/d,c'=c/d$ and $d=gcd(b,c)$ and we put $b'=1,c'=0$ if $c=0$.
In case $c=b-1$ there is a unique genus one fibration with class
$\epsilon_1=(1,0)$.

In case $b>2,c=1$, the non-trivial automorphism of $Y_{b,c}$
interchanges the two elliptic fibrations.
More generally, if $Y_{b,c}$ has two elliptic fibrations
it is interesting to compare
the associated relative Jacobians.

\subsection{Proposition.}
Let $b\in \ZZ_{>0}$ and $0\leq c<b-1$ and assume that $gcd(b,2c)=1$.
Then the relative Jacobians of the two elliptic fibrations on
$Y_{b,c}$ are isomorphic.

In case $(b,c)=(2,0)$ the relative Jacobians are not isomorphic in general.

\ts
If $gcd(b,2c)=1$, the discriminant group of the transcendental lattice
$T\cong \Gamma_{b,c}$ of $Y_{b,c}$ is cyclic of
order $b^2$.
Hence it has a unique subgroup $H$ of order $b$, which is
maximal isotropic.
Thus $\Gamma_{b,c}$ has a unique even, unimodular overlattice
and this overlattice must be $\Gamma_{1,0}\cong U^2\oplus E_8(-1)^2$.
Hence there is a unique Hodge structure on $\Lambda_{K3}$
such that $T$ is a sublattice of index $b$ in the transcendental lattice
of this Hodge structure. As this Hodge structure determines a unique
elliptic fibration, this must be
the relative Jacobian fibration (for any genus one fibration on $X$).

The surface $Y_{2,0}$ is a 2:1 cover of $\PP^1\times\PP^1$ branched over
a curve of bidegree $(4,4)$ (cf.\ section \ref{2,0}).
If the Jacobian fibrations associated to the two elliptic pencils
were isomorphic,
there would be an automorphism of $\PP^1$ carrying the discriminant locus
of one fibration into the one of the other fibration, which is easily seen
not to be true in general.
(In this case there are two maximal isotropic subgroups
in the discriminant group $\ZZ/2\ZZ\times\ZZ/2\ZZ$
of $T$ and the induced Hodge structures on the
two overlattices are not necessarily Hodge isometric).
\qed

\subsection{Example}\label{example}
Let $b>2$ be a prime number.
Given $\alpha\in Br(X)$ of order $b$, the subgroup $<\alpha>$ has $b-1$
elements of order $b$.
Since $\alpha$ and $-\alpha$ define the same genus one fibration,
there are $(b-1)/2$ genus one fibrations
$Y\rightarrow \PP^1$ with
$T_Y\cong T_{<\alpha>}$ (so these surfaces are Fourier-Mukai partners)
and with associated relative Jacobian $X/\PP^1$.

We consider the case $T_{<\alpha>}\cong \Gamma_{b,1}$
with $b\equiv 1\rmod 4$, so that $-1$ is a square mod
$b$, in some more detail.
Proposition \ref{Lbc} shows that the N\'eron-Severi group isomorphic
to $\Lambda_{b,-1}$, $\Lambda_{b,1}$ or $\Lambda_{b,c^2}\cong
\Lambda_{b,c^{-2}}$ with $c^2\not\equiv c^{-2}\rmod b$,
hence we get $1+1+(b-5)/4$ isomorphism classes.
If $NS(Y)\cong\Lambda_{b,-1}$ there is a unique genus one fibration on $Y$.
If $NS(Y)\cong\Lambda_{b,1}$ there are two genus one fibrations on $Y$,
the automorphism of the surface $Y$ permutes them, so they are isomorphic.
The other surfaces have two genus one fibrations.
These fibrations are not isomorphic in general
but both have relative Jacobian fibration $X/\PP^1$,
so we get $1+1+2((b-5)/4)=(b-1)/2$ genus one fibrations associated
to $<\alpha>$ as expected.
We conclude that there are $1+1+(b-5)/4=(b+3)/4$ isomorphism
classes of K3 surfaces $Y$ with $T_Y\cong T_{<\alpha>}$.

\section{Some projective models of K3 surfaces with genus one fibrations}
\label{exaY}

\subsection{} The well-developed theory of linear series on
K3 surfaces (cf.\ \cite{SD}) allows one to construct projective
models for
a K3 surface $Y_{b,c}$
with N\'eron-Severi group isomorphic to $\Lambda_{b,c}$
as in \ref{ybc}.
For small $b$ we recover some well-known examples.

\subsection{The case $c=b-1$}\label{cb-1}
When $c=b-1$ and $b\geq 3$,
the class $h_0=(0,1)\in C^+_{b,c}$, with $h_0^2=2c$, defines a map:
$$
\phi_{h_0}:Y_{b,c}\longrightarrow \PP^{c+1}.
$$
The unique genus one fibration $\epsilon_1$
has class $(1,0)$.
As $h_0\epsilon_1=b$ the curves from this genus one pencil
are mapped to degree $b=c+1$ elliptic normal curves. Each curve spans
a $\PP^{c}\subset \PP^{c+1}$, so it spans a hyperplane section.
The residual intersection is $h_0-\epsilon_1=N$, which is the $(-2)$-curve.
As $h_0N=2c-b=c+1$, $\phi_{h_0}(N)$ is a rational normal curve in $\PP^c$
and it meets each elliptic curve from $\epsilon_1$ in $N\epsilon_1=c+2$
points. See also the examples \ref{2,1} and \ref{3,2}.

\subsection{The case $c=0$}\label{c=0}
When $c=0$, the class $h=(1,1)$ maps $Y=Y_{b,0}$ to $\PP^{b+1}$.
As $h=\epsilon_1+\epsilon_2$, the sum of the genus one fibrations,
if we fix a curve $E_1\in |\epsilon_1|$, the residual intersection
of any hyperplane containing $E_1$ is a curve from $|\epsilon_2|$.
As $h\epsilon_1=b$, the curve $E_1$ is an elliptic curve of degree $b$
which spans a $\PP^{b-1}$ in $\PP^{b+1}$,
so we get a pencil of hyperplanes ($\PP^{b}$'s) whose residual intersections
give the system $|\epsilon_2|$. In particular, fixing
independent sections $s_1,s_2\in H^0(Y,\epsilon_1)$,
$t_1,t_2\in H^0(Y,\epsilon_2)$ we get the sections
$z_{ij}=s_it_j\in H^0(Y,h)$ which satisfy $z_{11}z_{22}=z_{12}z_{21}$,
hence the image of $Y$ under $\phi_h$
is contained in a rank 4 quadric whose rulings
cut out the linear systems $|\epsilon_1|$ and $|\epsilon_2|$.
See the examples \ref{2,0}, \ref{3,0}.

\subsection{$(b,c)=(1,0)$}\label{1,0}
There is a unique $-2$-curve $N$ on $Y_{1,0}$
with class $(-1,1)$ and there is a unique genus one
fibration $F$ with class $(1,0)$. Since $NF=1$, the curve $N$ is
a section of the fibration, so $Y_{1,0}$ has
an elliptic (Jacobian) fibration.

Note that if $D\in C^+_{1,0}$ is very ample, then we must have $DF\geq 3$
(else the fibers are not embedded) so the class of $D$ is $(a,b)$ with
$a>b\geq 3$ and thus $D^2\geq 12$.

\subsection{$(b,c)=(2,0)$}\label{2,0}
The two elliptic pencils give
a 2:1 map from $Y_{2,0}$ to $\PP^1\times \PP^1$ branched over a curve
of bidegree $(4,4)$.
The Segre map embeds $\PP^1\times\PP^1$ as a smooth quadric in $\PP^3$
(cf.\ \ref{c=0}).

\subsection{$(b,c)=(2,1)$}\label{2,1}
Let $D=\epsilon_1+N=(0,1)$, so $D^2=2$.
The 2:1 map $\phi_D:Y\rightarrow \PP^2$
contracts the $(-2)$-curve $N=(-1,1)$, the inverse images of the lines
through $\phi_D(N)$ are the union of curves in the genus one fibration
and $N$.
Thus $Y_{2,1}$ is obtained as the desingularization of the double cover of
$\PP^2$ branched along a sextic with one node. In particular, $Y_{2,1}$
is a double cover of the ruled surface $\FF_1$ branched along a smooth curve
curve of genus $9$.

\subsection{$(b,c)=(3,0)$}\label{3,0}
The K3 surface $Y_{3,0}$ is
a complete intersection of a rank 4 quadric and a cubic hypersurface in
$\PP^4$ (cf.\ \ref{c=0}).
The rank 4 quadric has two rulings, i.e.\
families of $\PP^2$'s parametrized by a $\PP^1$, each such $\PP^2$
intersects the surface along the curve where it intersects the cubic
hypersurface.

\subsection{$(b,c)=(3,1)$}\label{3,1}
A general K3 surface $Y$ of bidegree $(2,3)$ in $\PP^1\times \PP^2$
has $NS(Y)\cong \Lambda_{3,1}$.
Let $x=\{pt\}\times\PP^2$, $y=\PP^1\times\PP^1$ in $Pic(\PP^1\times \PP^2)$,
then $x^2=0,xy=\{pt\}\times\PP^1,y^2=\PP^1\times \{pt\}$
and $xy^2=1,y^3=0$. Restricting these
classes to the K3 surface one finds the lattice $\Lambda_{3,1}$
(for example, $y^2$ restricts to $y^2(2x+3y)=2$).
This surface has an obvious elliptic fibration (projection onto $\PP^1$)
and is also a double cover of $\PP^2$ (projection on $\PP^2$).
Writing the equation of the surface as
$$
Y:\quad s^2F+2stG+t^2H=0\qquad (\subset \PP^1\times\PP^2)
$$
with homogeneous
polynomials $F,G,H\in\CC[x_0,x_1,x_2]$ of degree 3, the branch curve of
the double cover is $G^2-FH$. In particular, this sextic has tangent
cubics (like $F=0$, in fact for any $s,t$, the cubic with equation as above
is tangent to the branch curve)
and these are the fibers of the elliptic fibration.
This system of tangent cubics defines an (even, effective)
theta characteristic on the branch curve.
The surface $Y$ has a second elliptic
fibration given by $-x+3y$ which is the image of the original fibration $x$
under the covering involution.
The Segre map $\PP^1\times \PP^2\rightarrow \PP^5$
embeds this surface in $\PP^5$.

\subsection{$(b,c)=(3,2)$}\label{3,2}
A general quartic surface $Y$ with a line in $\PP^3$  has
$NS(Y)\cong \Lambda_{3,2}$.
The unique
genus one fibration is given by the pencil of
planes on the line.
The general complete intersection of type (1,1) and (1,3) in
$\PP^1\times\PP^3$ is also of this type.
If $x=\{pt\}\times\PP^3$ and $y=\PP^1\times\PP^2$
then in degree 4 the only non-zero intersection is $xy^3=1$.
The class of the surface is $(x+y)(x+3y)=x^2+4xy+3y^2$. Hence
restricting $x,y$ to the surface one finds the lattice $\Lambda_{3,2}$.
The Segre map
$$
\phi_h:Y\longrightarrow \PP^1\times\PP^3\hookrightarrow \PP^6\subset\PP^7,
$$
is given by $h=x+y$.
The equation of type $(1,1)$ is a linear equation for the
image of $Y$ in the Segre embedding, hence the image lies in a hyperplane
$\PP^6$ in the Segre space $\PP^7$.

\subsection{$(b,c)=(4,2)$}
The general complete intersection $Y$ of bidegree
(1,2) and (1,2) in $\PP^1\times\PP^3$ has $NS(Y)\cong \Lambda_{4,2}$.
If $x=\{pt\}\times\PP^3$
and $y=\PP^1\times\PP^2$ then in degree 4 the only non-zero intersection
is $xy^3=1$. The class of the surface is $(x+2y)^2=x^2+4xy+4y^2$, hence
restricting $x,y$ to the surface and finds the lattice $\Lambda_{4,2}$.
The Segre map is $\phi_h:Y\rightarrow \PP^7$ with $h=x+y$.
Using the coordinates $(\lambda:\mu)$ on $\PP^1$, the two equations
for the surface are:
$$
\lambda Q_1-\mu Q_2,\qquad \lambda Q_3-\mu Q_4,
$$
for homogeneous degree two polynomials $Q_1,\ldots,Q_4$
in the coordinates of $\PP^3$. In
particular, the projection onto $\PP^3$ has equation $Q_1Q_4-Q_2Q_3=0$
which is a K3 surface containing the elliptic curves defined by
$Q_1=Q_3=0$ and $Q_1=Q_2=0$. These quartic genus one
curves lie in distinct pencils.

\section{Brauer groups of double covers of rational surfaces}\label{brdo}

\subsection{} In this section we prove a general theorem relating
the 2-torsion in $Br(X)$ of a double cover $X$ of a rational surface
to the 2-torsion of the Jacobian $J(C)$ of the branch curve $C$.
The theorem and its proof are due to the referee.
The theorem greatly generalizes
a similar theorem for elliptic fibrations in a previous version
which used results from \cite{KPS}.
Besides the Example \ref{exthmn}, we will apply this
theorem in \ref{brj} and \ref{h^2=2}.

\subsection{Theorem}\label{thmn}
Let $f:X\rightarrow Y$ be a double cover of a smooth complex
projective rational surface $Y$ branched along a smooth curve $C$.
Let $G\cong \ZZ/2\ZZ$ be the subgroup of $Aut(X)$ generated by the covering
involution $\sigma$. Assume that $\sigma$ acts trivially on $Pic(X)$
and that $Pic(X)$ is torsion-free.

\noindent
Then there is a natural inclusion
$$
J(C)_2\;\hookrightarrow\; Br(X)_2
$$
with quotient $(\ZZ/2\ZZ)^n$ and $n=2(1+b_2(Y)-b_0(C))-\rho$
where $\rho=\text{rk}NS(X)$ and $b_i$ is the $i$-th Betti number.

\ts
We use the spectral sequences in the category of $G$-sheaves on the $G$-variety
$X$ (\cite{Grot}, Chap.\ V):
$$
E^{p,q}_2=H^p(Y,R^qf_*^G{\cal F})\Longrightarrow H^n(G;X,{\cal F})
$$
$$
\pE^{p,q}_2=H^p(G,H^q(X,{\cal F}))\Longrightarrow H^n(G;X,{\cal F}),
$$
to determine $H^2(G;X,{\cal F})$ for ${\cal F}=\cO^*_X$.
We have $R^0f_*^G\cO^*_X=\cO_Y^*$ and for $q>0$
the sheaf $R^qf_*^G\cO^*_X$ on $Y$ has support
on the branch curve $C$ with stalks
$$
(R^qf_*^G\cO^*_X)_y=H^q(G,\cO^*_{X,y})
$$
where we identify $y\in C$ with its inverse image in $X$
(\cite{Grot}, Th\'eor\`eme 5.3.1).
As $Y$ is rational, we have $E^{2,0}_2=H^2(Y,\cO_Y^*)=0$, hence also
$E^{2,0}_\infty=0$.
Thus the filtration on $E^2=H^2(G;X,\cO^*_X)$ has only two steps
(cf.\ \cite{CE}, Chap.\ XV, Prop.\ 5.5) and is given by the exact sequence
$$
0\longrightarrow E_\infty^{1,1}\longrightarrow E^2\longrightarrow
E_\infty^{0,2}\longrightarrow 0.
$$
For a $G$-module $A$ we have $H^0(G,A)=A^G$, the submodule of invariants,
$H^{2q+1}(G,A)=ker(1+\sigma)/im(1-\sigma)$ and, for $q>0$,
$H^{2q}(G,A)=ker(1-\sigma)/im(1+\sigma)$.
Thus $H^2(G,\cO^*_{X,y})$ is the quotient of $\cO^*_{Y,y}$
by the subgroup generated by the
$\delta_1(h)=h\cdot (h\circ \sigma)$ for $h\in\cO^*_{X,y}$.
A germ $[g\in\cO^*_Y(U)]\in \cO^*_{Y,y}$, with $U$ simply connected,
has a square root on $U$. Let $h\in\cO^*_X(f^{-1}U)$ be the
the composition of such a square root with $f$.
Then $h=h\circ\sigma$ and $g=\delta_1(h)$
so $H^2(G,\cO^*_{X,y})$ is trivial, for all $y\in Y$.
Therefore $E^{0,2}_\infty=0$. Next we show that
$$
H^1(G,\cO^*_{X,y})\stackrel{\cong}{\longrightarrow} \mu_2:=\{1,-1\},
\qquad
g\longmapsto g(y).
$$
An element in $ker(1+\sigma)$ is a germ $g\in\cO^*_{X,y}$ satisfying
$g\cdot (g\circ\sigma)=1$. As $\sigma(y)=y$ this implies
$g(y)=\pm 1$. Moreover, if $g\in im(1-\sigma)$ then
$g=h\cdot (h\circ \sigma)^{-1}$ for some $h\in\cO^*_{X,y}$
and thus $g(y)=1$, so the map above is well defined.
It is obviously surjective (take $g=-1$). If $g$ is in the kernel,
so $g(y)=1$, let $h=1+g$, note that $h\in \cO^*_{X,y}$.
Using $g\cdot (g\circ\sigma)=1$, we have
$$
g\cdot(h\circ \sigma)= g\cdot(1+g\circ\sigma)
=g+g\cdot(g\circ\sigma)=1+g=h
$$
so $g=h\cdot (h\circ \sigma)^{-1}$,
and thus $g$ is trivial in $H^1(G,\cO^*_{X,y})$.
Therefore we get the
isomorphism (note that we followed the proof of `Hilbert 90' closely).
Hence $R^1f_*^G\cO^*_X$ is the constant sheaf $\mu_2$ on $C$ and
$E^{1,1}_2=H^1(C,\mu_2)=J(C)_2\cong (\ZZ/2\ZZ)^{b_1(C)}$.

To compute $E^{1,1}_\infty$ recall that (since $E^{p,q}_2=0$ when
$p$ or $q$ is negative):
$$
E^{1,1}_\infty=E^{1,1}_3=
ker(d^{1,1}_2:\,E^{1,1}_2\longrightarrow\, E^{3,0}_2).
$$
The exponential sequence
and the fact that $Y$ is a surface, so $H^3(Y,\cO_Y)=0$, show that
$E^{3,0}_2=H^3(Y,\cO_Y^*)\cong H^4(Y,\ZZ)=\ZZ$.
In particular, $E^{3,0}_2$ is torsion-free,
so the map $d^{1,1}_2$ must be trivial
and hence $E^{1,1}_\infty=E^{1,1}_2=J(C)_2$. Therefore the first spectral
sequence gives
$$
H^2(G;X,\cO^*_X)\cong J(C)_2.
$$

In the second spectral sequence we have
$$
\pE^{2,0}_2=H^2(G,H^0(X,\cO^*_X))=H^2(G,\CC^*)=0
$$
since $G$ acts trivally on $\CC^*$ whence $ker(1-\sigma)=\CC^*=im(1+\sigma)$.
Thus $\pE^{1,1}_\infty$ is a subgroup of $\pE^2$
with quotient $\pE^{0,2}_\infty$.
By assumption, $G$ acts trivially on $Pic(X)=H^1(X,\cO_X^*)$
and $Pic(X)$ is torsion-free. Therefore $\pE^{1,1}_2=H^1(G,Pic(X))=0$
because $\ker(1+\sigma)$ is trivial. Thus also $\pE^{1,1}_\infty$ is trivial
and $\pE^2\cong\pE^{0,2}_\infty$.
Next we have
$$
\pE^{0,2}_\infty=\pE^{0,2}_3=
ker(d^{0,2}_2:\pE_2^{0,2}\longrightarrow \pE_2^{2,1})=
ker(H^2(X,\cO_X^*)^G\longrightarrow H^2(G,Pic(X))).
$$
As $\pE^2=H^2(G;X,\cO^*_X)=J(C)_2$ is a 2-torsion group,
its image in $H^2(X,\cO_X^*)^G$  is also 2-torsion. This image is
also a subgroup of $H^2(X,\cO_X^*)$, so it is a subgroup of
$Br(X)_2$. This gives the desired inclusion
$J(C)_2\hookrightarrow Br(X)_2$.

Finally we observe that the Hurwitz formula for the cover $f$ gives:
$$
\chi_{top}(X)=2\chi_{top}(Y)-\chi_{top}(C),
$$
as $b_1,b_3$ of the surfaces involved are zero and as $C$ is a disjoint
union of Riemann surfaces, so $b_0(C)=b_2(C)$ is the number of components
of $C$, we get:
$$
 2+b_2(X)=2(2+b_2(Y))-(2b_0(C)-b_1(C))\quad\Longrightarrow\quad
b_2(X)=2(1+b_2(Y)-b_0(C))+b_1(C).
$$
As $Br(X)_2\cong (\ZZ/2\ZZ)^m$ with $m=b_2(X)-\rho$, one finds that
$Br(X)/J(C)_2\cong(\ZZ/2\ZZ)^n$ with $n=m-b_1(C)=2(1+b_2(Y)-b_0(C))-\rho$.
\qed

\subsection{Example}\label{exthmn}
Let $X=Y_{2,0}$ as in \ref{2,0}, in particular $X$ is a K3 surface
(so $b_2(X)=22$) and $\rho(X)=2$. Thus $Br(X)_2\cong (\ZZ/2\ZZ)^{20}$.
Then $X$ is a double cover
of $Y=\PP^1\times \PP^1$ branched over a smooth, irreducible curve $C$
of bidegree $(4,4)$, hence the genus of $C$ is $9$ and $J(C)_2\cong
(\ZZ/2\ZZ)^{18}$. Therefore we get an exact sequence
$$
 0\longrightarrow J(C)_2\longrightarrow Br(X)_2\longrightarrow (\ZZ/2\ZZ)^2
\longrightarrow 0.
$$
Note that $b_2(Y)=2$ and and $b_0(C)=1$.

In the example \ref{2,1} one has a K3 surface $X=Y_{2,1}$ which is a
double cover of the rational ruled surface $Y=\FF_1$ branched along a
smooth irreducible curve $C$ of genus 9. Hence also in this case
one finds an exact sequence as above.

\section{Two torsion in the Brauer group of an elliptic surface}
\label{ttbr}

\subsection{Ruled surfaces}
We denote by $\bar{f}:\FF_e\rightarrow \PP^1$ the ruled surface
which is the compactification of the line bundle $\cO(e)$ on $\PP^1$
with a rational curve $C_\infty$ at infinity,
with self-intersection $C_\infty^2=-e$.
This surface can be obtained, for $e\geq 1$,
as the blow up of the cone in $\PP^{e+1}$
over the rational normal curve of degree $e$ in the vertex,
in case $e=0$ we have $\FF_0\cong\PP^1\times\PP^1$.

Let $F$ be a fiber of $\bar{f}$,
then
$$
Pic(\FF_{e})=\ZZ F\oplus \ZZ C_\infty,\qquad
F^2=0,\quad FC_\infty=1,\quad C_\infty^2=-e.
$$
The image in $\FF_{e}=\cO(e)\coprod C_\infty$
of a global section of $\cO(e)$ is a section $C$
of $\FF_e\rightarrow \PP^1$.
It is also a hyperplane section of the cone in $\PP^{e+1}$
which doesn't meet the vertex.
The classes of the canonical line bundle $K$ on $\FF_{e}$
and the class of a section $C$  of $\bar{f}$ are:
$$
K=K_{\FF_{e}}=-(e+2)F-2C_\infty,\qquad C=eF+C_\infty
$$
($K$ can easily be found using the adjunction formula and the fact that
$F$ and $C_\infty$ have genus zero, the class of $C$ can be found using
$CF=1$ and $CC_\infty=0$).

\subsection{Weierstrass models of elliptic surfaces}\label{wms}
Let $f:X\rightarrow \PP^1$ be a general elliptic fibration (with section).
Then $X$ is the double
cover of the ruled surface $\bar{f}:\FF_{2d}\rightarrow\PP^1$
for  some $d\in\ZZ_{\geq 0}$ (\cite{CD} $\S$5.5, \cite{FM} 1.4.1, \cite{Mi}).
The branch locus of $\pi:X\rightarrow\FF_{2d}$ consists of the section at
infinity $C_\infty$ of $\bar{f}$
and a curve $\cC$ such that $\pi:\cC\rightarrow \PP^1$ is
a 3:1 covering of $\PP^1$.
The topological Euler characteristic of $X$ is given by the Hurwitz formula:
$$
\chi(X)=2\chi(\FF_{2n})-\chi(\cC\coprod C_\infty)=8-(2-2g+2))=4+2g.
$$
Since $H^1(X,\ZZ)=H^3(X,\ZZ)=0$, the rank of $H^2(X,\ZZ)$ is
$2+2g$.

We assume from now on that $Pic(X)\cong\ZZ^2$,
so we get
$$
Br(X)\cong H^2(X,\QQ/\ZZ)/(Pic(X)\otimes(\QQ/\ZZ))\cong (\QQ/\ZZ)^{2g}
$$
and $\cC$ is smooth and irreducible.
From the intersection numbers $C_\infty \cC=0$, $F\cC=3$ one finds that
$\cC$ has class $3C$.
The canonical bundle $K_\cC$ on $\cC$ is the restriction of $K+\cC$ to $\cC$,
that is of $(4d-2)F+C_\infty$ to $\cC$, this gives:
$$
\cC=3C=6dF+3C_\infty,\qquad g(\cC)=6d-2.
$$

\subsection{The theta characteristic on $\cC$}\label{defthc}
As the canonical bundle on $\cC$ is the restriction of
$K+\cC=(4d-2)F+C_\infty$ to $\cC$ and $C_\infty$ and $\cC$ are disjoint,
the canonical bundle of $\cC$ is also the restriction of $(4d-2)F$ to $\cC$.
Let $\theta$ be the restriction of the divisor class $(2d-1)F$ to $\cC$:
$$
\theta=(2d-1)F_{|\cC},\qquad{\rm then}\qquad 2\theta=K_\cC
$$
hence $\theta$ is a theta characteristic on $\cC$.

\subsection{Lemma}\label{lemthc}
The theta characteristic $\theta$ on $\cC$ is even
and $\dim H^0(\cC,\theta)=2d$.

\ts
Let $D=\cC\cap F$, a divisor of degree three on $\cC$ with $h^0(D)=2$,
let $1,x$ be a basis of $H^0(D)$ (so the map $\bar{f}:\cC\rightarrow\PP^1$
is given by these two functions).
As $\theta=(2d-1)D$, we get $1,x,\ldots,x^{2d-1}\in H^0(\theta)$.

For $y\in H^0(\theta)$, we get
$y,xy,\ldots x^{2d-1}y\in H^0(2\theta)=H^0(K_\cC)$,
this space also contains $1,x,\ldots,x^{4d-2}$.
Hence if $y$ is not a linear combination of powers of $x$,
then
$\dim H^0(K_\cC)\geq 2d+4d-1=6d-1$, but
$\dim H^0(K_\cC)=g(\cC)=6d-2$, so we get a contradiction.
\qed

\subsection{Remark}
The theta characteristic $\theta$ on $\cC$
gives a natural bijection between the points of
order two in $Pic(\cC)$ and the theta characteristics of $\cC$ given by
$\alpha\mapsto \alpha+\theta$.

\subsection{Theorem}\label{brj}
Let $f:X\rightarrow \PP^1$ be an an elliptic fibration
with Weierstrass model
$\pi:X\rightarrow \FF_{2d}$,
branch curve $\cC\coprod C_\infty$,
and $Pic(X)\cong\ZZ^2$.
In particular, $\cC$ is a smooth curve of genus $g=6d-2$.

Then there is a natural isomorphism:
$$
Br(X)_2\cong J(\cC)_2.
$$

\ts This follows immediately from Theorem \ref{thmn} applied to
$\pi:X\rightarrow \FF_{2d}$. Note that the branch
curve $C$ of $\pi$ has two components, so $b_0(C)=2$,
and that $Pic(X)$ is generated by the class of a section and a fiber, so
the covering group $G$ of $\pi$ acts trivially on $Pic(X)$.
\qed

\section{Conic bundles over Jacobian fibrations}\label{cobu}

\subsection{}
In this section we consider a genus one fibration
$g:Y\rightarrow \PP^1$ with a double section
and we recover the associated Jacobian fibration
$f:X\rightarrow \PP^1$. For this we only use the classical construction
of the Jacobian of a double cover
of $\PP^1$ branched over four points, due to Hermite (cf.\ \cite{Weil},
\cite{Mu}), i.e.\ we consider the generic fibers of $g$ and $f$ over
the field $k(\PP^1)$. This construction also gives us a conic bundle
over $X$ which corresponds to $\alpha\in Br(X)$ defining $Y$.

\subsection{The Jacobian fibration}
Let $g:Y\rightarrow\PP^1$ be a a genus one fibration defined by an
equation
$$
Y:\qquad w^2=a_0v^4+4a_1v^3+6a_2v^2+4a_3v+a_4
$$

The Jacobian fibration of $g$ was found by Hermite (cf.\ \cite{Weil}).
We will follow \cite{Mu}, 1.3b and 11.3c.
The Jacobian fibration of $f:Y\rightarrow\PP^1$ as above is
the elliptic surface $f:X\rightarrow \PP^1$ given by
$$
X:\qquad y^2=4x^3-g_2x-g_3,
$$
with
$$
\begin{array}{rcl}
g_2&=&a_0a_4-4a_1a_3+3a_2^2,\\
&&\phantom{x}\\
g_3&=&a_0a_2a_4-a_0a_3^2-a_1^2a_4+2a_1a_2a_3-a_2^3.
\end{array}
$$
In case $Y$ is a K3 surface,
the $a_i$ are homogeneous polynomials and the
$g_2$ and $g_3$ are
homogeneous of degree $8$ and $12$ respectively
(cf.\ \ref{doqu} and \ref{nodsex}),
so the K3 surface $X$ maps onto $\FF_{4}$.

We will assume that the curve
$$
\cC:\qquad 4x^3-g_2x-g_3=0\qquad(\subset\FF_{2d})
$$
in the ruled surface $\FF_{2d}$ is smooth.
The surface $X$ is a double cover of $\FF_{2d}$, the
branch curve in $\FF_{2d}$ is the union of $\cC$ and the `section at
infinity' $C_\infty$.

\subsection{}\label{matM}
The explicit expressions for $g_2$ and $g_3$ imply that the equation
of $\cC$ is the determinant of a symmetric matrix $M$:
$$
4x^3-g_2x-g_3=\det(M),\qquad
M=\left
(\begin{array}{ccc}
a_0&a_1&a_2+2x\\a_1&a_2-x&a_3\\a_2+2x&a_3&a_4
\end{array}
\right).
$$

Thus $M$ defines a conic bundle on $\FF_{2d}$ which degenerates over the
curve $\cC$. In particular, it defines an \'etale 2:1 cover of $\cC$
(the conics split in two lines over points in $\cC$), thus we get
a point $\alpha$ of order two in $Pic(\cC)$ (cf.\ \cite{Bepr}).
From sections \ref{b=2}, \ref{2,0}, \ref{2,1} it follows that
the theta characteristic $\alpha + \theta$,
where $\theta$ is the even theta characteristic on
$\cC$ defined in \ref{defthc}, is odd in the case of a double $\FF_1$
and even in the
case of a double quadric.

Alternatively, one could try to define $\alpha$
as the cokernel of a map of locally free sheaves $M: \cG\rightarrow \cF$
(as in the case of symmetric matrices with entries in $\CC[x_0,x_1,x_2]$
where one considers bundles over $\PP^2$, cf.\ \cite{Be}, \cite{BCZ}).

\subsection{Double quadrics}\label{doqu} (See \ref{2,0}.)
Let $G\in\CC[s,t,u,v]$ be bihomogeneous of degree $4$ in $s,t$
and in $u,v$, so
$G(\lambda s,\lambda t,\mu u,\mu v)=\lambda^4\mu^4G(s,t,u,v)$.
Then $G$ defines a curve $\cR$ of bidegree $(4,4)$ and genus $9$
in $\PP^1\times\PP^1\cong\FF_0$.
Let $\bar{f}:\PP^1\times\PP^1\rightarrow \PP^1$ be the projection to the
first factor. The double cover of $\PP^1\times\PP^1$
branched over $\cR$ is a K3 surface $Y$ with a genus one
fibration $f:Y\rightarrow \PP^1$,
the elliptic curves are the inverse images of the fibers of
$\bar{f}$.
This fibration can be defined by an equation of the form:
$$
w^2=a_0v^4+a_1v^3+a_2v^2+a_3v+a_4,\qquad a_i\in\CC[s,t],\qquad deg(a_i)=4,
$$
where $G(s,t,u,v)=a_0(s,t)v^4+\ldots+a_4(s,t)u^4$. Note that each
$a_i$ is homogeneous of degree $4$.

\subsection{Nodal sextics}\label{nodsex} (See \ref{2,1}.)
Let $F\in\CC[x_0,x_1,x_2]$ be a homogeneous polynomial of degree six which
defines a plane sextic $F=0$ with a node in $(0:0:1)$.
Blowing up $\PP^2$ in this point, one obtains the ruled surface
$\FF_1\subset\PP^1\times\PP^2$, the map $\bar{f}:\FF_1\rightarrow \PP^1$
is the projection on the first factor.
The double cover of $\FF_1$ branched over the strict transform $\cR$
of the sextic (which has genus $10-1=9$) is a K3 surface $Y$ with a genus one
fibration $f:Y\rightarrow \PP^1$,
the elliptic curves are the inverse images of the fibers of
$\bar{f}$.
This fibration is defined by an equation of the form:
$$
w^2=a_0v^4+a_1v^3+a_2v^2+a_3v+a_4,\qquad a_i\in\CC[s,t],\qquad deg(a_i)=i+2,
$$
where $F(x_0,x_1,x_2)=a_0(x_0,x_1)x_2^4+\ldots+a_4(x_0,x_1)$. Note that
$a_i$ is homogeneous of degree $i+2$.

\subsection{The trigonal construction}
In the two cases above, the branch curve $\cR\subset\FF_e$, $e=0,1$,
of $Y\rightarrow \FF_e$ is a curve of genus $9$.
The fibration
$\FF_e\rightarrow \PP^1$ induces a 4:1 map $\cR\rightarrow \PP^1$, i.e.\
$\cR$ is a tetragonal curve.
The genus one fibration on $Y$ corresponds to a point
of order two $\alpha\in Br(X)_2\cong J(\cC)_2$, where $X\rightarrow \PP^1$
is the relative Jacobian of $Y/\PP^1$ and $\cC$ is the branch curve of
$X\rightarrow \FF_4$. The fibration $\FF_4\rightarrow \PP^1$ induces
a 3:1 map $\cC\rightarrow \PP^1$, so we get a trigonal curve $\cC$
with a point of order two $\alpha\in J(\cC)_2$.
The pairs $(\cC,\alpha)$ correspond to \'etale 2:1 covers
$\tilde{\cC}_\alpha\rightarrow \cC$.

The relation between $(\cC,\alpha)$ and $\cR$
is given by Recillas' trigonal construction (\cite{R}, \cite{D} \S 2.4)
which gives a natural bijection between \'etale
double covers of trigonal curves
and tetragonal curves. In fact the trigonal construction associates to
a 4:1 cover of $\PP^1$,
that is to an extension of degree four of $k(\PP^1)$,
the cubic extension of $k(\PP^1)$ given by the Lagrange resolvent and
its natural quadratic extension.
The explicit formula for the Lagrange resolvent
(cf.\ \cite{Go}, II, $\S$174)
shows that the cubic extension corresponds to the
3:1 cover $\cC\rightarrow \PP^1$.

\subsection{An Azumaya algebra}
The symmetric matrix $M$ in \ref{matM} defines a conic bundle
on $\FF_{2d}$ which ramifies over $\cC$,
so the fibre over each point of $\cC$ is the union of two lines.
This provides an unramified double cover of $\cC$, or, equivalently,
a point of order two on $J(\cC)$ (cf.\ \cite{Bepr}).
The pull-back along $\pi:X\rightarrow \FF_{2d}$
gives an unramified conic bundle on $\pi^{-1}(\FF_{2d}-\cC)$.
Due to the branching of order two along $\cC$, this
conic bundle extends to all of $X$
(locally near a point of $X$ on the ramification locus, the
conic bundle is defined by a homogeneous polynomial of the
form $X^2+Y^2+t^2Z^2$, where $t=0$ defines the ramification curve.
Changing coordinates $Z:=tZ$ shows that the conic bundle does not ramify in
codimension one, hence does not ramify at all).

The even Clifford bundle associated to the (unramified)
conic bundle on $X$ is an Azumaya algebra, actually it
is a bundle of quaternion algebras (cf.\ \cite{Pu} $\S$2).

The (generic fiber of the) conic bundle on $X$ is determined by a
symmetric $3\times 3$ matrix with coefficients in the field of functions
$k(X)$ of $X$. Diagonalizing this matrix we obtain a matrix
$diag(f,g,h)$. The associated even Clifford algebra is the quaternion
algebra with symbol $(-fg,-fh)$, i.e.\ it is the rank 4 algebra over
$k(X)$ with generators ${\bf i},\bf{j}$ and relations
${\bf i}^2=-fg$,
${\bf j}^2=-fh$, ${\bf i}{\bf j}=-{\bf j}\bf{i}$.

Conversely, to check that such an element of $Br(k(X))$
defines an element of $Br(X)$ one can use
the tame symbol (cf.\ \cite{W}).
Wittenberg's explicit example \cite{W} of elements of order two in
the Brauer group of an elliptic fibration is related to a `degenerate'
case (the curve $\cC$ has three irreducible components)
of the construction considered in this section.

\section{Remarks on general K3 surfaces.}\label{genk3}

\subsection{Two torsion in $Br(X)$.}\label{h2=2}
Let $X$ be a K3 surface with $NS(X)=\ZZ h$ and $h^2=2d>0$.
Then $Br(X)_2\cong (\ZZ/2\ZZ)^{21}$.
As in the case of elliptic fibrations
(cf.\ Proposition \ref{lemed}),
there are various types of elements $\alpha\in Br(X)_2$
corresponding to the isomorphism classes of the sublattices
of index two
$$
\Gamma_\alpha=\ker(\alpha:T_X\longrightarrow \ZZ/2\ZZ),
\qquad\alpha\in Br(X)_2,\;\alpha\neq 0.
$$

To classify these sublattices we observe that the inclusion
$NS(X)=\ZZ h\hookrightarrow \Lambda_{K3}$ is unique up to isometry.
Hence we may assume:
$$
NS(X)=\ZZ h\cong \ZZ(1,d)\hookrightarrow U
\hookrightarrow U\oplus \Lambda'\cong H^2(X,\ZZ),\qquad
{\rm with}\quad
\Lambda'= U^2\oplus E_8(-1)^2.
$$
Let $v=(1,-d)\in U$, so $v^2=-2d$. Then the transcendental lattice is:
$$
T:=T_X\cong \ZZ v\oplus \Lambda'\cong \langle-2d\rangle
\oplus U^2\oplus E_8(-1)^2.
$$
Since $\Lambda'$ is a unimodular (and hence self-dual) lattice,
any $\alpha\in Br'(X)_2$ can be written as:
$$
\alpha:T\longrightarrow\ZZ/2\ZZ,\qquad
nv+\lambda'\longmapsto a_\alpha n+\langle\lambda_\alpha,\lambda'\rangle
\;\rmod 2
$$
for a unique $a_\alpha\in\{0,1\}$ and a $\lambda_\alpha\in\Lambda'$
whose class in $\Lambda'/2\Lambda'$ is uniquely determined
by $\alpha$.

To classify the index two sublattices of $T_X$
we only need to determine their discriminant groups and
discriminant forms (\cite{N}, 1.13.3).

\subsection{Proposition}\label{propd}
Let $X$ be a K3 surface with $NS(X)=\ZZ h$ and $h^2=2d$.
Let $\alpha\in Br'(X)_2=Hom(T,\ZZ/2\ZZ)$, $\alpha\neq 0$,
be defined by $a_\alpha\in\{0,1\}$ and
$\lambda_\alpha\in\Lambda'$ as in \ref{h2=2}.
Let
$\Gamma_\alpha=\ker(\alpha)$ be the index two sublattice of $T$ defined
by $\alpha$. Then:
\begin{itemize}
\item[  i)]
$\Gamma_\alpha^*/\Gamma_\alpha\cong
\ZZ/2d\ZZ\oplus\ZZ/2\ZZ\oplus\ZZ/2\ZZ$ if $a_\alpha=0$,
there are $2^{20}-1$ such lattices.
\item[ ii)]
$\Gamma_\alpha^*/\Gamma_\alpha\cong \ZZ/8d\ZZ$ if $a_\alpha=1$,
there are $2^{20}$ such lattices.
\begin{itemize}
\item[ a)]
In case $d\equiv 0\rmod 2$, these $2^{20}$ lattices are isomorphic
to each other.
\item[ b)]
In case $d\equiv 1\rmod 2$, there are two isomorphism classes of
such lattices. One class, the even class, has $2^9(2^{10}+1)$
elements and is characterized by
$\mbox{$\frac{1}{2}$}\langle\lambda_\alpha,\lambda_\alpha\rangle
\equiv 0\rmod 2$, the other (odd) class has $2^9(2^{10}-1)$
elements and is characterized by
$\mbox{$\frac{1}{2}$}\langle\lambda_\alpha,\lambda_\alpha\rangle
\equiv 1\rmod 2$.
\end{itemize}
\end{itemize}

\ts
Since $\alpha\neq 0$, the discriminant of $\Gamma_\alpha$ is
$2^2{\rm disc}(T)=8d$, and this is the order of
$\Gamma_\alpha^*/\Gamma_\alpha$. Let $a=a_\alpha$.

In case $a=0$, we have an orthogonal direct sum decomposition
$\Gamma_\alpha=\ZZ v\oplus \Gamma'_\alpha$
where $\Gamma'_\alpha=\ker(\alpha)\cap\Lambda'$.
Hence $\Gamma_\alpha^*/\Gamma_\alpha$ is the direct sum of the
discriminant group of $\ZZ v$ (which is $\ZZ/2d\ZZ$) and the
discriminant group of $\Gamma'_\alpha$, which has order 4.
It is easy to see that
$\lambda_\alpha/2\in \Gamma_\alpha^*/\Gamma_\alpha$
and that it has order two.
Next choose a $\mu\in\Lambda'$ with
$\langle\lambda_\alpha,\mu\rangle=1$, then also
$\mu\in \Gamma_\alpha^*/\Gamma_\alpha$ has order two. This proves
the case $a=0$.

In case $a=1$, the element $u:=(-v,2d\lambda_\alpha)/4d\in T\otimes\QQ$
is in the lattice $\Gamma_\alpha^*$ since for
$(nv,\lambda')\in\Gamma_\alpha$:
$$
\langle u,(nv,\lambda')\rangle=
\mbox{$\frac{1}{4d}$}(2dn+2d\langle\lambda_\alpha,\lambda'\rangle)=
\mbox{$\frac{1}{2}$}(n+\langle\lambda_\alpha,\lambda'\rangle)
= \mbox{$\frac{1}{2}$}\alpha(nv,\lambda') \in \ZZ.
$$
Note that
$4du=(-v,2d\lambda_\alpha)
\not\in\Gamma_\alpha$ but $8du\in
\Gamma_\alpha$. Hence $u$ has order $8d$ in the discriminant
group and thus it is a generator.
The value of the discriminant form on this
generator is:
$$
q(u)=
q\bigl(\frac{(-v,2d\lambda_\alpha)}{4d}\bigr)=
\frac{-2d+4d^2\langle\lambda_\alpha,\lambda_\alpha\rangle}{16d^2}=
\frac{-1+2d\langle\lambda_\alpha,\lambda_\alpha\rangle}{8d}
\rmod 2\ZZ.
$$
To find the isomorphism classes, we need to
see if there exists an integer $x$ such that
$$
x^2\cdot \frac{-1}{8d}\equiv \frac{-1+4d}{8d}\quad\rmod 2\ZZ.
$$
Multiplying by $8d$ we get the equivalent condition:
$$
x^2\equiv 1-4d\quad \rmod 16d.
$$
If $d$ is odd,
this gives a contradiction when considered modulo $8$.
In case $d$ is even, we write $d=2^re$ with $e$ odd and $r\geq 1$.
Then the equation becomes:
$$
x^2\equiv 1-2^{r+2}e\quad \rmod 2^{r+4}e.
$$
Since $e$ is odd, the equation is equivalent to $1-2^{r+2}e$ being a square
modulo $2^{r+4}$ (which is true: any integer $\equiv 1\rmod 8$ is a square
modulo any power of two) and modulo $e$ (which is trivial:
$x^2\equiv 1-2^{r+2}e\equiv 1\rmod e$ is satisfied by $x=1$).
\qed

\subsection{} In contrast to the case when $X$ has an elliptic fibration,
not all of the sublattices $\Gamma_\alpha$ allow a primitive embedding
into $\Lambda_{K3}$. In case $\Gamma_\alpha$ has a primitive embedding,
this embedding is not necessarily unique (\cite{O}, \cite{St}).

\subsection{Corollary} The lattice $\Gamma_\alpha$ has a primitive
embedding into $\Lambda_{K3}$ if and only if $a_\alpha=1$ and,
in case $d\equiv 1\rmod 2$,
$\mbox{$\frac{1}{2}$}\langle\lambda_\alpha,\lambda_\alpha\rangle
\equiv 0\rmod 2$.

\ts
If $\Gamma_\alpha$ has a primitive embedding $i$ into $\Lambda_{K3}$,
then we must have $i(\Gamma_\alpha)^\perp\cong \ZZ f$ for some $f$
with $f^2=8d$
for reasons of rank and discriminant. Thus such an $i$ exists if
and only if the discriminant form of $\Gamma_\alpha$ is minus the one of
$\ZZ f$. The computation done in the proof of Proposition \ref{propd}
gives the result. \qed

\subsection{The case $h^2=2$.}\label{h^2=2}
If $h^2=2$, the map $\phi_h:X\rightarrow \PP^2$
is a double cover branched over a
smooth sextic curve $C\subset\PP^2$, so the genus of $C$ is ten.
From Theorem \ref{thmn} it follows that there is an exact sequence:
$$
0\longrightarrow J(C)_2\stackrel{}{\longrightarrow}
Br(X)_2\longrightarrow \ZZ/2\ZZ\longrightarrow 0.
$$
It was shown in \cite{F}, Theorem 2.4, that there is an exact sequence:
$$
0\longrightarrow Br(\PP^2-C)_2\stackrel{\phi_h^*}{\longrightarrow}
Br(X)_2\longrightarrow \ZZ/2\ZZ\longrightarrow 0.
$$
In \ref{remaining} we will see that there is natural way to associate
a conic bundle on $\PP^2-C$ to an element in $J(C)_2$, so we get
a natural isomorphism $J(C)_2\cong Br(\PP^2-C)_2$.

From the results we recall below it also follows that the group $Br(X)_2$
is naturally isomorphic to the subgroup of $Pic(X)$ generated
by the theta characteristics, modulo the subgroup generated by the canonical
class (note that if $\theta$ and $\theta'$ are theta characteristics, then
$\theta+\theta'=K+L$ for some $L\in J(C)_2$ because
$2(\theta+\theta')=2K$).

\subsection{The case $h^2=2$ and $\alpha$ is even}
If $d=1$, $a_\alpha=1$ and
$\mbox{$\frac{1}{2}$}\langle\lambda_\alpha,\lambda_\alpha\rangle
\equiv 0\rmod 2$, the lattice $\Gamma_\alpha$ has a unique primitive
embedding into $\Lambda_{K3}$ and thus defines a unique K3 surface
$Y_\alpha$. It is easy to see that $Y_\alpha$ has a polarization
of degree $8$ and that
 $Y_\alpha$ has an embedding into $\PP^5$ where it is the
complete intersection of three quadrics (cf.\ \cite{Mub}, \cite{Kh}).
The discriminant curve $C$ of this
net of quadrics is the plane sextic over which $X\rightarrow \PP^2$
branches. The net of quadrics defines, and is defined by,
an even theta characteristic on $C$ (cf.\ \cite{Be}).
Thus one has a canonical bijection between the even theta characteristics
on the genus ten curve $C$ and the `even' sublattices
$\Gamma_\alpha\subset T_X$ (see also \cite{La}, \cite{OG} for this
sublattice of index two). In \cite{Bh} it is shown that $Y_\alpha$ is the
moduli space of certain orthogonal bundles on $X$: this is a kind of inverse
to Mukai's construction of $X$ as a moduli space of bundles on $Y_\alpha$.

\subsection{The case $h^2=2$ and $\alpha$ is odd}
If $d=1$, $a_\alpha=1$ and
$\mbox{$\frac{1}{2}$}\langle\lambda_\alpha,\lambda_\alpha\rangle
\equiv 1\rmod 2$, the lattice $\Gamma_\alpha$ does not have a primitive
embedding into the K3 lattice. However, it is isometric to a primitive
sublattice of $H^4(Z_\alpha,\ZZ)$ where $Z_\alpha\;(\subset\PP^5)$
is a cubic fourfold which contains a plane $P$. In fact,
$$
\Gamma_\alpha(-1) \cong \langle h^2,P\rangle^\perp
\qquad(\subset H^4(Z_\alpha,\ZZ)),
$$
where $h\in H^2(Z_\alpha,\ZZ)$ is the hyperplane class, cf.\ \cite{V},
and the quadratic form on $\Gamma_\alpha(-1)$ is the opposite to the one
on $\Gamma_\alpha$.
The projection of the fourfold $Z_\alpha$ from the plane $P$ defines
a quadric bundle structure on $Z_\alpha$. The rulings of the quadrics
define a double cover of $\PP^2$, which is $X$.
The quadric bundle defines, and is defined by, an odd
theta characteristic on $C$ (cf.\ \cite{V}, \cite{La}, \cite{Be})
and in this way one obtains a
canonical bijection between the odd theta characteristics
on the genus ten curve $C$ and the `odd' sublattices
$\Gamma_\alpha\subset T_X$. For more interesting relations between
K3 surfaces and cubic fourfolds see \cite{hassett}.

\subsection{The case $h^2=2$ and $a_\alpha=0$}\label{remaining}
In case $d=1$ and $a_\alpha=0$, we have $\alpha\in Br(\PP^2-C)_2$.
The conic bundle $P_\alpha$ over
$X$ defined by $\alpha$ is obtained by pull-back from a conic
bundle $Q_\alpha$ on $\PP^2$ which ramifies along $C$.
Thus $Q_\alpha$ defines an \'etale 2:1 cover of $C$, that is,\
an invertible sheaf $L_\alpha\in J(C)_2-\{0\}$.

Conversely, the resolution of the sheaf $i_*L_\alpha(1)$ on
$\PP^2$, here $i:C\hookrightarrow \PP^2$ is the inclusion,
gives a symmetric $3\times 3$ matrix $M_\alpha$ with entries
which are homogeneous of degree two
and whose determinant vanishes on $C$. Thus $M_\alpha$ defines the
conic bundle $Q_\alpha$ (cf.\ \cite{BCZ}). In this way one obtains
canonical bijections between the points of order two in $J(C)$,
the points of order two in $Br(\PP^2-C)$
and the sublattices
$\Gamma_\alpha\subset T_X$ with $a_\alpha=0$.

For $x,y\in\CC^3$ we thus have a quadratic form ${}^tyM_\alpha(x)y$,
which is a polynomial of bidegree (2,2) in $x$ and $y$. Let
$W_\alpha\subset \PP^2\times \PP^2$ be the threefold defined by
${}^tyM_\alpha(x)y=0$. The projections $\pi_x$ and $\pi_y$
define two conic bundle structures on $W_\alpha$.
In particular, we get two sextics, $C=C_x$ and $C_y$
as branch curves and thus we get two K3 surfaces,
branched over these curves, $X$ and $X_\alpha$.

On the other hand, the lattice
$\Gamma_\alpha\cong <-2>\oplus U(2)\oplus U\oplus E_8(-1)^2$
has exactly two embeddings
as index two sublattice in the lattice
$T$.
In fact, these embeddings correspond to isotropic subgroups
of the discriminant form of $\Gamma_\alpha$.
As $\Gamma_\alpha$ has discriminant form
$(\ZZ/2\ZZ)^3\rightarrow \QQ/2\ZZ$,
$(\bar{x},\bar{y},\bar{z})\mapsto \mbox{$\frac{1}{2}$}x^2+yz\rmod 2\ZZ$,
the isotropic subgroups are those generated by $(0,1,0)$ and $(0,0,1)$.

\

\end{document}